\documentclass[12pt,leqno]{article}
\title{Updating An Upper Bound Of Erik Westzynthius} 

\author{Gerhard R. Paseman} 
\begin{document}

\maketitle
\newcommand{\gt}{>}
\newcommand{\lt}{<}
\newcommand{\sigin}{{\sigma}^{-1}(n)}
\newcommand{\piin}{{\pi}^{-1}(n)}
\newcommand{\soik}{\sum_{1 \leq i}^k}
\newcommand{\poik}{\prod_{1 \leq i}^k}
\newcommand{\eqdf}{=_{\textrm{\small{ def }}}}
\newcommand{\myfloor}[1]{\lfloor #1 \rfloor}
\newcommand{\myceil}[1]{\lceil #1 \rceil}

\abstract{Inspired by a paper of Erik Westzynthius,
 we build on work of Harlan Stevens and Hans-Joachim 
Kanold.  Let $k \gt 2$ be the number of distinct prime divisors 
of a positive integer $n$.  In 1977, Stevens used Bonferroni
inequalities to get an explicit upper bound on Jacobsthal's 
function $g(n)$, which is related to the size of largest interval of 
consecutive integers none of which are coprime to $n$. 
Letting $u(k)$ be the base $2$ $\log$ of this bound, Stevens
showed $u(k)$ is $O((\log k)^2)$, improving upon Kanold's
exponent $O(\sqrt{k})$.  We use elementary methods similar to those of Stevens 
to get $u(k)$ is $O(\log k(\log\log k))$ in one form and 
$O(\sigin\log k)$ in another form.  We also show how 
these bounds can be improved for small $k$.}

\section{Overview}

Erik Westzynthius provided a ground-breaking result in \cite{W} on 
prime gaps, showing that for any constant $D \gt 0$ there were 
infinitely many primes $p_n$ so that $p_{n+1} \gt p_n + D\log(p_n)$.
In the same paper, he also provided an upper bound for what was to
be later called $g(P_k)$.  This quantity measures how many
consecutive integers we can find having a "small" $( \leq p_k)$
prime factor.

Ernst Jacobsthal in \cite{J} defined and showed
$g(n) \leq k2^k + 2^k - k$, where $k$ is the 
number of distinct prime factors of $n$.  Improved explicit bounds
were later given by Kanold \cite{K} ($2^k$ for all $k$ and $2^{\sqrt{k}}$ 
for $k \geq e^{50}$)
and Stevens \cite{S} $(2k^{2 + 2e\log k})$, and some additional
but less explicit results given by Paul Erd\H{o}s, Henryk Iwaniec, and others.

We show several explicit bounds for $g(n)$, some depending on the 
quantities $\sigin$ and $\piin$, including
\begin{displaymath} 
g(n) \lt \frac{\sum_{i=1}^s {k \choose i}}{\piin - (\sigin)^{s+1} /(s+1)!},
\end{displaymath}
with $s$ an odd integer bounded by $(1$ plus a constant times $\sigin)$,
which gives $g(n) \lt k^{B + C\log\log k}$ for $k>2$ with explicit
constants $B$ and $C$, both $\lt 3.9$.

We define $\sigin$ and $\piin$ in the next section, and also list
some results which apply when $\sigin$ is small $(\lt 1 + 1/2q_1)$. 
The following section recalls work of Jacobsthal and
Westzynthius, and shows how Stevens's bound can be tightened with
a little effort. 

With elementary means, we also show the bound above that uses $\sigin$ and
$\piin$, and follow that section with some supplementary results as
well as suggestions for further research.  The remainder of this 
article contains some history, a recommended reading list, and an
Appendix as well as acknowledgments and a list of citations.

\section{Definitions and simple bounds}

We use $\omega(n)$ for the number of distinct (positive)
prime factors of the positive integer $n$, and declare
$k=\omega(n)$.  We also require $k \gt 0$, so $n \gt 1$.
We list the prime factors $q_i$ of $n$ in increasing order:
$q_1 \lt q_2 \lt \ldots \lt q_k$. 
We recall the $k$th primorial $P_k$ as $P_k=\poik p_i.$
 
\textbf{Definition}: We define (here $^{-1}$ is
part of a label, \textbf{not} an exponent)
$$\sigin = \soik 1/q_i \textrm{ , and } 
\piin = \poik (1 - 1/q_i).$$
Note that $n\piin = \phi(n)$, Euler's function for counting
positive integers coprime to and less than $n$.

\textbf{Definition}: After Jacobsthal, for $n \gt 0$ 
define  $g(n)$ to be the smallest positive integer 
$m$ such that for any integer $a$, the set of $m$ consecutive 
integers $\{a+1,\ldots, a+m\}$ has an integer $a+j$, where 
$1 \leq j \leq m$, such that $\gcd(a+j,n)=1$.  If $n \gt 1$,
define $L(n)$ as the largest integer $l$ so that there is an 
interval of $l$ consecutive integers $\{b+1,\ldots, b+l\}$ such
that each $b+j$ satisfies $\gcd(b+j,n) \gt 1$.

It is straightforward to show $g(n)=L(n)+1$ for all 
$n \gt 1$.  Also $g(n) = \max_i (c_{i+1} - c_i)$ where
the $c_i$ denote all integers coprime to (\textbf{totatives} of)
 $n$ in increasing
order.  Consequently, $g(n)$ depends only on the set
of distinct prime factors of $n$; we will use $n$
squarefree at times in this article.

One has $g(p_i)=2$, and $g(n) \lt n$ for $n \gt 2$.  If $n$
has "large" prime factors ($q_i \gt k$ for all $i$), then
$g(n)=L(n)+1$ is "small".  Specifically,

\textbf{Proposition}(Jacobsthal):    If $q_1 \gt k$, then $L(n)=k$.

Proof Sketch:  Any integer interval of length $k+1$ has at 
most one multiple of $q_i$, giving at most $k$ integers in 
that interval having a prime factor in common with $n$.  
So at least one integer in that interval is not a 
multiple of any $q_i$. So $L(n) \lt k+1$.

Conversely, for any permutation $\tau$ of the $k$ indices, 
the Chinese Remainder Theorem gives an integer $b_\tau$
 such that $b_\tau + \tau(i) = 0 \bmod{q_i}$ for
$1 \leq i \leq k$, so $L(n) \geq k$.  \textbf{End} of Proof Sketch.

We could have chosen $\tau$ to be the identity permutation, but
we want to point out that this proof gives $k!$ many intervals 
in $(0,n)$ which achieve the $L(n)$ bound for this kind of $n$.
This proof also shows $g(n) \gt k$ for any $n \gt 1$.  

The next bound is inspired by Kanold's 1967 paper.  We would like to
know if it appears explicitly in the literature.

\textbf{Proposition}(Kanold-\textbf{P}.): Let $0 \lt r \lt 1$ be a real number 
and $n$ such that $r + \sigin \lt 1$.  Then 
$L(n) \lt k/r$, and $g(n) \leq \myceil{k/r}$. 

Proof Sketch:  For an interval containing $L$ consecutive 
integers, at most $\myceil{L/q_i} \leq 1 + (L-1)/q_i$ of 
them are multiples of $q_i$.  Since $\sigin \lt 1-r$, 
 the count of numbers not coprime to $n$ is at most 
$(L-1)\sigin + k \lt (L-1)(1-r) +k$.  Now whenever $L \geq k/r$,  
$(L-1)(1-r) + k \leq L -1 + r$, so this count is less than $L$. 
One gets $g(n) \leq \myceil{k/r}$, 
giving $L(n) \lt k/r$.  \textbf{End} of Proof Sketch.

This gives a weaker bound than Jacobsthal's proposition, but
it applies to more cases, even when $q_1$ is about $2\sqrt{k}$
and $r \gt 1/q_1$.  It suggests the following

\textbf{Variation}: Let $k \gt 1$ and $\sigin \lt 1 + \frac{1}{2q_1}$. Then 
$$L(n) \lt \frac{q_1}{q_1-1}(2k -1 - \sigin)2q_1 \textrm{, so } 
g(n) \lt 4kq_1\frac{q_1}{q_1-1}.$$

Proof Sketch:  The estimate above for multiples of $q_i$,
$1 \lt i \leq k$, among the $L$ numbers is refined by subtracting that  
portion which are also multiples of $q_1$; we
underestimate it by 
$L/q_iq_1 - 1 \lt \myfloor{\myceil{L/q_i}/q_1}$.
Thus $(L-1)\sigin + k - \sum_{1 \lt i}^k (-1 + L/(q_iq_1))$ is the
refined
upper bound for the non-totative count.  Rewriting, we look for $L$ so that 
\begin{displaymath}
L[\sigin(1- 1/q_1) + 1/(q_1)^2] + 2k-1 - \sigin \lt L,
\end{displaymath}
which holds if and only if 
$$L(1 - 1/q_1)(1 + 1/q_1 - \sigin) \gt (2k - 1 -\sigin).$$
The assumption $(1 +1/q_1 - \sigin) \gt 1/2q_1$ yields that
if $L/2q_1 \gt (2k - 1 - \sigin)(q_1/(q_1 -1))$, then
$L$ is greater than our overcount; as above, we get
$L(n) \lt 2q_1 (2k - 1 - \sigin)(q_1/(q_1 - 1))$ 
and the weaker $g(n) \lt 4kq_1(q_1/(q_1 -1))$.  \textbf{End} of Proof Sketch.

We now present one bound of general character, and one which asymptotically
improves upon Kanold's smaller bound.  We use the  

\textbf{Fact}: Let $b, d$ and $f$ be integers with $\gcd(f,d)=1$.  
There is an integer $z$ so that for all integers 
$t, \gcd(b + tf, d) = \gcd(zb + t, d)$.

Pick $z$ so that $zf=1 \bmod{d}$. $\gcd(b + tf, d) = \gcd(z(b + tf), d) = \gcd(zb + t, d)$,
giving the Fact above and the application below that coprimality
in sequences of consecutive integers behaves the same way in
certain arithmetic progressions. 

Assume $k \gt 2, n$ squarefree and $d \mid n$ with
$d$ neither $1$ nor $n$.  Set $f=n/d$.  $1+tf$ is coprime to $f$, and
the Fact above shows that $1+tf$ always has at least one
out of $g(d)$ consecutive members coprime to $d$.  Any interval of $(g(d)f)$-many
integers thus has one coprime to $n$ of the form $1+tf$. 
Considering $c+tf$ for all totatives $c$ of $f$, we get the

\textbf{Observation}: $g(fd) \leq g(d)f - f + g(f).$

Proof Sketch: When $c \neq 1$, 
there is at least one number of the form $c+tf$ coprime to
$n=fd$ in the interval $(1+a, 1+a + g(d)f)$, where $a=t_0f$ and we
assume that if $1+tf$ is coprime to $n$, then it is outside this interval.
There are at least $(\phi(f) - 1)$-many of these numbers; try placing them
inside the interval leaving a large gap.  If you can get a gap
of at least $g(d)f - f$, then the numbers must be 
in the subintervals $(1+a, 1+a+f)$ and 
$(1+a+(g(d)-1)f, 1+a+g(d)f)$.  If the largest in the first
interval is $c+a$, then the smallest in the second interval is
at most $c'+a+(g(d)-1)f$, where $c'$ is the next largest totative
to $f$ after $c$.  But $c' - c \leq g(f)$, giving the bound.
\textbf{End} of Proof Sketch.

We use weaker forms of the Observation and Variation to improve on
Kanold's smaller bound.

Let $n$ be squarefree.  Consider the tail of $\sigin$,
that is, find $q_l$ smallest such that 
$\sum_{l \leq i}^k 1/q_i \lt 1 + 1/2q_l$.  Let
$d = \prod_{l \leq i}^k q_i$.  From the Variation we have 
$g(d) \lt 4(k-l+1)({q_l}^2/(q_l - 1))$.  However, when $k$ is large, 
$q_l$ is small enough that $n/d$ is smaller 
than $2^{\sqrt{k}}$.

\textbf{Improvement}:  For $k$ sufficiently large,
$n/d \lt 2^{\frac{3}{2}k^{0.45}}$, giving
$g(n) \lt 2^{3k^{0.45}/2}\frac{4(k-l+1){q_l}^2}{q_l-1}$.
As $k$ gets large, $\log g(n) \lt k^{(\epsilon +1/e)}$ 
eventually for any fixed $\epsilon \gt 0$.

Proof Idea:  Let $m$ satisfy $p_m=q_{l-1}$. We defer to the Appendix showing 
that $p_m \lt k^{0.45}$ when $k \geq e^{9.5}$.
As $n/d \leq P_m$, and for an explicit
positive $c \lt 0.02$, $P_m \lt e^{p_m(1+c)} \lt 2^{3p_m/2}$; 
the Observation implies the weaker form 
$g(n) \lt (n/d)g(d) \lt 2^{3k^{0.45}/2}g(d)$.

About $\log g(n)$, we use that 
$\sum_{l \leq j}^k 1/q_j \leq \sum_{1\leq i}^{k-l+1} 1/p_{m+i}$, and 
then truncate to the first $u \leq k-l+1$ terms while still having
$\sum_{1 \leq i}^u 1/p_{m+i} \gt 1 - 1/2p_m$.  We use an 
approximation of Mertens for this last sum to get 
$\log(\log(p_{m+u})/\log p_m ) \gt 1 - O(1/\log p_m)$. Picking
a convenient positive $\epsilon_1$ and $\epsilon$, as 
$m$ goes large, $p_{m+u} \gt  p_m^{e - \epsilon_1}$ eventually, 
and $\log$ of the bound is majorized by $p_m(1+c')$ for some
positive small $c'$ 
and is eventually dominated by $k^{1/e + \epsilon}$.
The Appendix provides more detail.  
\textbf{End} of Proof Idea.
 
We take credit for this form of presentation, but are 
influenced by Kanold's paper;  his proof for $2^{\sqrt{k}}$
has many of the ideas above, and we wonder if perhaps
he did discover it.

We end this section with lower bounds:
Westzynthius shows an easily demonstrated lower bound 
$2p_{k-1} \leq g(P_k)$  for $k \gt 1$, and then shows a better
lower bound (for $k$ sufficiently large after choosing $\epsilon \gt 0$) of 
$e^{\gamma}(1- \epsilon) p_k \frac{\log\log p_k}{\log\log\log p_k}$. 

\section{Improving estimates of Stevens}

Both Westzynthius and Jacobsthal use a simple sieve argument to
establish an upper bound on $g(n)$, using $n$ squarefree.  
Since Stevens uses part of 
this, we show the argument here.  

Recap: We use inclusion-exclusion to count integers in the open interval
$L=(a, a+x)$ that are coprime to squarefree $n$,
 for real numbers $a$ and 
$x \gt 0$.  $I_0$ counts the totatives of $n$ in $L$, 
 and for $t \gt 0$ and $t \mid n,$ $I_t$ 
counts multiples of $t$ in $L$. 
Then $I_0 = \sum_{t\mid n} (-1)^{\omega(t)}I_t.$
As in the previous section we replace the count
$I_t$, this time by $x/t + E(t)$.  $E(t)$ is an error term which depends 
on both $t$ and $a$ actually, but always $|E(t)| \leq 1,$
and will be removed below.
$$I_0 = \sum_{t \mid n} (x/t + E(t))(-1)^{\omega(t)} \geq 
x \sum_{t \mid n} (-1)^{\omega(t)}/t - \sum_{t \mid n} 1 .$$ 
Now rewrite the sum of $(-1)^{\omega(t)}/t$ as a product 
$\prod (1 - 1/q_i)$, and note the second sum is $2^k$. We get
$$I_0 \geq x\piin -2^k .$$

Now this relation above is essentially independent of $a$.  If 
we pick $x$ large enough, then $I_0 \gt 0$ and $g(n)$
will be at most $x$.  So choose $x = 2^k/\piin + \epsilon$.
\textbf{End} of Recap.
 
The above shows $g(n) \leq 2^k/\piin$.  One would like to improve
on this since Kanold has an asymptotically better result.

In our view Stevens has two main ideas in his 1977 paper.
The first is to observe that the 
sum produced from inclusion-exclusion can be truncated
a la Bonferroni to a smaller sum to give fewer terms to 
approximate.  The second is that the denominator $T_s$ can 
be written as $\piin - T_s'$, and that $T_s'$ can be easier to handle.  
We adapt his proof slightly, and then we tighten
up the estimates he uses to improve the exponent. 

Adaptation: We use $x, I_0, t,$ and $I_t$ where Stevens used $Q, L, B$ and $N(\ldots).$
We use integral $x$, following Stevens. We assume $k \gt 4$
to make sure some of his estimates apply. He reorganizes
the sum by the number of factors in $t$ and uses a result of Landau
for the first idea.  His display (3) in our notation says: for any
odd value of $s$,
$$I_0 \geq \sum_{0\leq i}^s (-1)^i \sum_{t \mid n, \omega(t)=i} I_t .$$

\newcommand{\myts}{\sum_{0 \leq i}^s (-1)^s \sum_{t \mid n, \omega(t)=i} \frac{1}{t}}  
\newcommand{\mysb}{\sum_{1 \leq i}^s {k \choose i}}
\newcommand{\dfeq}{_{\textrm{\small{def}}}=}

Using the estimate $| I_t - x/t | \leq 1$ (except for $t=1$
when $I_1 = x$), we write what Stevens has in his (4) and (5) as
$$I_0 \geq  x \myts - \mysb \eqdf xT_s - SB.$$
(We've written $SB$ for $\mysb$, and $T_s$ for $\myts$.)  Defining $T_s'$ by
the relation $\piin - T_s' = T_s$, Stevens notes (using an approximation
of $k^s$ for $SB$ instead of $SB$ directly) that if $s$ is chosen so that
$\piin \gt T_s'$ (so that $T_s \gt 0$), and if $x \gt SB/(\piin - T_s')$,
then $I_0 \gt 0$ and so $g(n) \leq x$.  He and we now look for a
suitable $s$.

$T_s'$ can be written as a sum over $i$ of sums of terms
$t$ with $\omega(t)=i$, just like $T_s$, but with $s \lt i \leq k$.
Before doing this, Stevens observes: for $k>2$
(we insert $h(k)$, a putative upper bound for $\sigin$)
$$r! \sum_{t \mid n, \omega(t)=r} 1/t \lt (\sigin)^r \lt h(k)^r .$$
Stevens uses $\log k$ for $h(k)$; later we will use $\log\log p_k + 1/2$.  Then
$$T_s' = \sum_{s \lt r}^k (-1)^r \sum_{ t \mid n,\omega(t)=r} 1/t
\lt  \sum_{s \lt r}^k \sigin^r / r! \lt \sum_{s \lt r}^\infty h(k)^r/r!$$
which follows by dropping the $(-1)^r$ and by extending the sum past $r=k$.
Taylor's theorem with remainder on $e^{h(k)}$ then yields
$$T_s' \lt \sum_{s < r}^\infty h(k)^r/r! \leq e^{h(k)} h(k)^{s+1}/(s+1)! .$$

Stevens bounds things further by asking
$s+1 \geq 2eh(k)$ and using $(s+1)! \gt ((s+1)/e)^{s+1}$.  Then
$$T_s' \lt e^{h(k)} h(k)^{s+1}/(s+1)! 
\lt e^{h(k)} (eh(k)/s+1)^{s+1} \leq e^{h(k)}/2^{s+1}.$$

He also under-estimates $\piin$ by $1/k$ for $k \gt 4$, 
where we will use $\beta/\log k$.  He then has 
$$\piin - T_s' \gt 1/k - k/2^{s+1} \gt 1/k - k^{1 - 2e\log 2} \gt 0,$$
since $s+1$ is an even integer greater than $2e \log k$.  
So $I_0 \gt 0$ (and thus $g(n) \leq x$) if $s$ is an odd 
integer greater than $2e\log k - 1$ and 
$x \gt 2k^{s+1} \gt k^s / ( 1/k - k/2^{s+1}) \gt SB/(\piin - T_s')$.
Stevens replaces $s+1$ with $2e\log k + 2$ 
to ensure the bound holds for all $k$.  \textbf{End} of Adaptation.

We repeat the above, using $\log\log p_k + 1/2$ for $h(k)$
 and $\beta/\log p_k$ for $\piin$. (The Appendix discusses the validity of
these choices.)
  Again asking for odd $s$ with
 $s+1 \geq 2eh(k),$ then  $T_s' \leq e^{h(k)}/2^{s+1}$, and
\begin{eqnarray*}
\piin - T_s' & \gt & \beta/\log p_k - (\log p_k)e^{1/2}/2^{s+1} \\
       & \gt & \beta/\log p_k - (\log p_k)e^{1/2}/2^{2e(\log\log p_k +1/2)} \\
       & =   & \beta/\log p_k - (\log p_k)(e/2^{2e})^{1/2} /(\log p_k)^{2e\log 2} \\
       & \gt & (1/\log p_k)[\beta - (\log p_k)^{2(1 - e\log2)}/3] \gt 0 .
\end{eqnarray*}
In the last line, we use that $(e/2^{2e})^{1/2} \lt 1/3$,
that we can pick $1/1.78 \gt \beta \gt 1/3$, that $\log p_k \gt 2$ because $k > 4$, 
and that $e\log 2 \gt 1$.  

It should be clear that $\piin - T_s' \gt 1/4\log p_k$ by choosing $\beta \gt 1/3$
and whenever $(\log p_k)^{2e\log 2 -2} > 4$ which holds for $k > 4$.  Using such an estimate
we have whenever $s+1$ is even and $ \geq 2e(\log\log p_k +1/2)$, $\piin -T_s' \gt 1/4\log p_k$ leading to
$$\textbf{Theorem: } g(n) \leq (4\log p_k) \sum_{1\leq i \leq s} {k \choose i} 
\lt (4\log p_k) k^{1 + 2e(1/2 + \log\log p_k)}.$$

We could tweak the choice of $s$ slightly to get a smaller exponent, as well as use a better
approximation for the sum of binomial coefficients.  In the next section, we will find a
bound which depends directly on $\sigin$ which not only does both, but gives a tighter bound
in general.

\section{$\sigin$ and $\piin$}

We modify Stevens's argument with a better upper bound
for the numerator and express $T_s'$ as an alternating
 and eventually decreasing sum, allowing us a smaller $s$. 

Note that $T_s'$ is an alternating sum and that one has
$$\sigin \sum_{t \mid n, \omega(t)=j} 1/t \gt (j+1) \sum_{t \mid n, \omega(t)=j+1} 1/t,$$
so that when $s > \sigin$, one can bound $T_s'$ by $\sigin^{s+1}/(s+1)!.$  Now instead
of Taylor's theorem and $h(k)$, we use (see Appendix)  
 $e \lt (1/\piin)^{1/\sigin} \leq 4$ to show any real number $r \geq 4\sigin$ gives
 $\piin - T_{\myceil{r} -1}' \gt \piin - \sigin^{\myceil{r}}/(\myceil{r})! \gt 0$:
\begin{eqnarray*}
 &  & e \lt 4^{3/4},\textrm{ so } e(1/\piin)^{1/4\sigin} \lt e4^{1/4}< 4\\
 & \textrm{so} & \sigin \lt \piin^{1/4\sigin}4\sigin/e \leq \piin^{1/r}(r/e) \\
 & \textrm{so} & \sigin^{\myceil{r}} \lt \piin(\myceil{r}/e)^{\myceil{r}} \lt \piin(\myceil{r})!.
\end{eqnarray*}

We now claim

\textbf{Theorem:} Let $s$ be the smallest odd integer with $s+1 \geq 4\sigin.$ For $k \gt 2,$
$$\frac{\sum_{1 \leq i \leq s} {k \choose i}}{\piin - \sigin^{s+1}/(s+1)!} \gt g(n).$$

We collapse the summands in the numerator slightly, increasing the total by 1, and
as $(s+1)! \gt \sqrt{2\pi(s+1)}((s+1)/e)^{s+1}$, one sees the denominator
is larger than $(\sqrt{2\pi(s+1)} - 1)\sigin^{s+1}/(s+1)!$, so we can write
a weaker upper bound as a \textbf{corollary}:
$$\frac{(s+1)!\sum_{0\leq 2j\lt s} { k+1 \choose s - 2j }}{(\sqrt{2\pi(s+1)} - 1)\sigin^{s+1}}
\gt g(n) .$$

This may seem intimidating, but when we take into account that for $k>6$, $\sigin$ is at most 
$1/2 + \log\log k ( 1 + \log 2/\log k)$
it is then seen that the dominant term in the sum is ${k+1 \choose s}$ and the
expression is $O(((k+1)/\sigin)^{3 + 4\log\log k + \epsilon})$ when $\sigin \geq 1$.
(The portion that is $\frac{s+1}{(\sqrt{2\pi (s+1)} -1)\sigin}$ is less than $1$ for
large enough $\sigin$; for $\sigin$ near or smaller than $1$ we have the more elementary
bounds.)

This argument only needs $K$ such that $(s+1) \gt K\sigin$ and also that
$e/K \lt \piin^{1/K\sigin}$.  This holds for $K \gt 3.89$, and when
$\sigin \gt 1$, one can lower $K$ from 4 to $3.81$.  However, even for large values of $\sigin$,
the argument still expects $(K/e)^K \gt e$, which means $K$ can't be shown
smaller than $3.59$ with this method.

As a rough comparison, Jacobsthal's bound is larger than Kanold's bound for all $k$.
Stevens's bound is smaller than $2^k$ for $k>300$, and is smaller than Kanold's better
bound for $k > 5000000$ .  Our Improvement is smaller than Kanold's for $k \gt e^{11}$,
and the exponent $s \leq 1 + 4\sigin$ above is smaller than that of Stevens for $k \gt 2$.

\section{Some History}

The recap is our interpretation of Westzynthius's upper bound argument published in 1931, 
generalized to arbitrary $n$ with $\omega(n)=k$ instead of $P_k$, which 
Westzynthius did not publish as far as we know.  In a footnote
Westzynthius did hint at sieving with just odd numbers, 
and we considered extending that argument with thinner sets.  This led us to asking the 
question \cite{M} on MathOverflow in 2010. 

Correspondence on MathOverflow led us eventually to Thomas Hagedorn's paper 
\cite{H} and Jacobsthal's function.
Jacobsthal in \cite{J}  uses a slightly different argument, and (with the
notation of this article) instead of 
using $2^k / \piin$ he bounds $\piin$ by $1/(k+1)$ and gives a bound of $(k+1)(2^k-1)$ on $L(n)$.
Hagedorn's paper quoted the bounds of Kanold and Stevens, and after studying those papers we
adapted Stevens's argument and posted the results on MathOverflow in 2011, as well as producing a
private manuscript with small circulation.  

Since then we have 
accumulated and posted other accessible results, and arranged some of them for this article.
The Observation represents a small improvement on Kanold's result which involves 
$g(fd) \leq g(d)f + 1- \phi(f)$
instead of $g(fd) \leq g(d)f + g(f)-f$; the two are the same for $f$ prime.
The Improvement is intended to
show not just improved asymptotic results but also that Kanold's bound holds for 
$k$ smaller than $e^{50}$. Indeed the name is earned once $k \geq e^{11}$.
We admit the work is in showing the bounds hold for small $k$, which makes the Improvement not as
elementary or accessible as we hope.

Hagedorn also mentions work of Erd\H{o}s, Iwaniec, and others.  Erd\H{o}s shows
for any given positive real $\epsilon$ that $| 1 - g(n)\piin/k | \gt \epsilon$ only for
$n$ in a set of zero density. Erd\H{o}s also comments that Brun's method can yield
a constant $c$ such that $g(n)$ is $O(k^c)$, but we have not found a version of this that is
both explicit and accessible.  Iwaniec shows the existence of a constant $C$ independent of $n$
such that there are at least $k^2$ many totatives of $n$ in an interval of size $C(k^2\log k)/\piin$,
which implies
$g(n)$ is $O((k\log k)^2)$; again we do not know what $C$ is.

\section{Further research and reading}

We intended this article to give simpler, more accessible, and more explicit proofs of upper
bounds on Jacobsthal's function.  We are optimistic about improving upon the
 results of Erd\H{o}s and Iwaniec.  In particular, we think there is more to the Observation:
we hope to achieve a subquadratic in $k$ upper bound using this direction by noting how large
intervals of numbers with  factor common to $n$ are distributed in $(0, 2n)$.  At present,
the difference in (base 2) exponents between $O(\log k\log\log k)$ and $O(\log k)$ is substantial.

Except for the bound depending on $\sigin$, all of these bounds are also bounds on
Jacobsthal's $C(k)$, given by $C(k) =\textrm{max}_{k=\omega(n)} L(n)$. It was shown recently 
\cite{HS} that
the conjecture $C(k) + 1 =g(P_k)$ holds for $1 \leq k \leq 23$ and fails at $k=24$.

Note that the bound involving $\sigin$ can represent a substantial improvement even
if $n$ cannot be factored; for those $n$ which do not have small factors, $\sigin$
can be substantially smaller as can $k$, even for numbers near $10^{10^{100}}$.  Of course,
when $k$ or a partial factorization of $n$ are better known, better bounds on $g(n)$
become available.

We are interested in tweaking the Variation to handle more squarefree $n$ 
by sieving out small prime factors.  Our beginning efforts have not yielded
much improvement on bounds obtained by the Observation.  It seems better
estimates on the number of totatives in an interval of arbitrary length are
needed to carry out an argument like that in the Variation.

Some questions of interest:

1) Pick a small odd prime $p$ and odd $n$ with $q_1 \gt p$.  We know $g(pn)/g(n) \lt p$: can we
get anything sharper?  In particular, what are those integers $n$ such that 
$g(3n) \gt 2g(n)?$  Such that $g(5n) \gt 2g(n)?$

2) Let $a(n)$ be the smallest positive integer such that $\gcd(n, a(n) + i) \gt 1$ for 
$0 \lt i \lt g(n)$.  One can show $a(n) \lt n/2$; how much can this be improved?  If $b(n)$
is the number of such longest intervals of nontotatives of $n$ in $(1,n)$, can we hope
for $a(n)b(n) \lt n$?

3) How close are two such intervals?  If one hopes for a subquadratic (in $k$) bound on
$g(n)$, this will be an important bit of information.  Even in the case $q_1 \gt k$,
it should be related to how primes are distributed, which suggests that some interesting
perspective is needed. 

4) Not much asymptotic improvement should be expected from these arguments in the case
that $\sigin$ is large, say $\sigin \gt 2$.  However, that is where the difficult cases
are, and the quantity $\cal{T}$  $= (1/\piin)^{1/\sigin}$ is expected to decrease as $\sigin$ increases.
How does $\cal{T}$ behave with $\sigin$, and can one use this in 
bounding $g(n)$?

5) Even the simple estimates with small $\sigin$ have some slop, primarily in overestimating
multiples of $q_i$ with
$L/q_i$.  Often this results in an estimate about twice as large as needed.  Can something
be said about this "noise" vector $L/q_i - I_{q_i}$ and what approaches avoid the error
introduced by this? 

We recommend the bibliography and also the following reading list, 
which provides additional information related to Jacobsthal's function and
applications.

Ernst Jacobsthal, \"{U}ber Sequenzen ganzer Zahlen, von denen keine zu $n$ teilerfremd ist. I-III.
\textit{Det Kongelige Norske Videnskabers Selskabs Forhandlinger} Bd 33 1960, Nr. 24,
Trondheim I Kommisjon Hos F. Bruns Bokhandel 1961, pp. 117-124,125-131,132-139.
(Also see IV and V published in a later edition.)
 
Hans-Joachim Kanold, Neuere Untersuchungen \"{u}ber die Jacobsthal- Funktion $g(n)$.
\textit{Monatshefte Math.} 84, 1977, pp. 109-124.

Paul Erd\H{o}s, On the integers relatively prime to $n$
and on a number-theoretic function considered by Jacobsthal.
\textit{Math. Scand.} 11 (1962) pp. 163-170.

R.C. Vaughan, On the order of magnitude of Jacobsthal's function.
\textit{Proc. Edinburgh Math. Soc.} 20 (1976-77) pp. 329-331.

Henryk Iwaniec, On the problem Of Jacobsthal.
\textit{Demonstratio Mathematica} v XI no. 1 1978 pp.225-231.

MathOverflow Questions (Number refers to URL, so 37679
expands to http://mathoverflow.net/questions/37679) \newline
88323 Analogues of Jacobsthal's Function. \newline
56099 Lower bound of the number of relatively primes (each other)
in an interval. \newline
68351 Least Prime Factor in a sequence of 2n consecutive integers. 

Shallit, J., State Complexity and Jacobsthal's Function.
\textit{CIAA '00 Revised Papers from the 5th International Conference on 
Implementation and Application of Automata} 2000, pp. 272-278

Schlage-Puchta, J.C., On Triangular Billiards. \newline
http://arxiv.org/abs/1105.1629 

\section{Acknowledgments}

We acknowledge Will Jagy, Aaron Meyerowitz, and Thomas Hagedorn for their support and 
assistance.  We appreciate the stimulating environment provided by the MathOverflow forum
and its community, and are thankful for its providing a repository of these results.
We especially appreciate Will and Aaron for their MathOverflow contributions, and are 
also grateful for discussions with 
users Larry Freeman (who asked a version of question 2) above) and asterios gantzounis.

\section{Appendix}

We resolve some details on assertions made in earlier sections: that
$(1/\piin)^{1/\sigin} \in (e,4]$, on bounding 
$\sigin$ from above by $1/2 + \log\log p_k$ and $\piin$ from below by
$1/3\log p_k$, and on showing the Improvement holds for $k \gt e^{9.5}$.

In getting a bound depending on $\sigin$, we used the assumption that
$e \lt (1/\piin)^{1/\sigin} \leq 4$ for $n \gt 1$.  We proved this
along with related results in a private manuscript \cite{P}. The 
proof was based on observing that
$(-\log \piin)/\sigin$ was a mediant sum of values of the form $p_i\log(p_i/(p_i-1))$
 ($\frac{a}{b}$ and
$\frac{c}{d}$ give a mediant of $\frac{a+c}{b+d}$), so that
$(q_k/(q_k-1))^{q_k} \lt (1/\piin)^{1/\sigin} \lt (q_1/(q_1-1))^{q_1}$.
Also, if $\sigin \geq 1$, then
$(1/\piin)^{1/\sigin}\lt 3.6$, so one can improve the constant $C$ in
$g(n) \leq Ak^{B + C\log\log k}$ to a value less than $3.81$.  If we did not care
about the advantage given by using $\sigin$, we could use a general bound of
$g(n) \lt k^{3 + 3.81\log\log k}$ for $k>2$, which can be verified by hand for small
values of $k$ and which would be weaker (and thus valid) than the Variation when $\sigin \leq 1$.

Letting $n=P_k$, Mertens determines $\sigin$ and $\piin$ with error by
$\sigin = \log\log p_k + B + E_1(k)$ and 
$\piin = e^{-(\gamma + \delta(k))}/\log p_k$ where
$B$ is a constant with value near $0.2615$ and
$E_1(k)$ and $\delta(k)$ are error terms in $O(1/\log p_k)$.  In using Stevens's argument
with tighter bounds, we used $1/2 + \log\log p_k$ as an upper bound for $\sigin$; calculations 
show that $\sigin \lt 0.41 + \log\log p_k$ for $7 \lt p_k \lt 10^8$, and Mertens estimates
of the error (or tighter estimates given by Rosser and Schoenfeld) show this holds for $k \gt 4$.
If we were concerned only with $k \gt 50$, we could replace $0.41$ by $0.28$, closer to the
value of $B$.  We could also use $3/4 + \log\log k$ as an upper bound for $k \gt 10$, and
replace $3/4$ by smaller numbers for $k$ sufficiently large.

Similarly $\piin \log(p_k)$ approaches $e^{-\gamma}$ which is near $0.5615$; computing
small examples shows $\piin \gt 1/3\log p_k$ for $1 \lt p_k \lt 10^8$; theory then gives the
weak inequality for all $k$.  We could use $1/2\log p_k$ for $k>8$ if we needed to improve the
multiplicative constant $4$; instead we chose to develop the estimate depending on $\sigin$
in the section following Stevens's argument.

The rest of this Appendix contains the verification of the claim
that the choice of $q_l$ for $n$ with $k \gt e^{9.5}$ satisfies $q_l \lt k^{0.45}$, and
remarks expanding on the proof idea of the Improvement.

We first work with sums of the form $\sum_{1 \leq j \leq u} 1/p_{m+j}$
which are within $1/2p_m$ to 1,
as they represent an extremal form with respect to the estimate. 
$\sum_{4 \leq j \leq 29} 1/p_i \lt .9$, so when $p_m=5$ we already have $u \gt p_m^2$ if we want
a sum close enough to 1. As one increases $m$ by 1, one
has to remove $1/p_{m+1}$ and "replace" it by more than $2p_m$ terms of size smaller than 
$1/4p_m^2$, so we already have $m/u \lt 1/m \log m$ for such sums.

As $m$ grows, $\log p_{m+u}$ tends to $e\log p_m$ and 
$(m+u) \log (m+u)$ approaches $p_m^e$, and thus $u/p_m^e$ approaches 1, yielding the
asymptotic (in $k$) result of $\log(g(n)) \lt k^{1/e + \epsilon}$.  Because of oscillations
around zero of the quantity $(\sum_{1 \leq j \leq u} 1/p_j -B - \log\log p_u)$, a proof
of $q_{l-1} \leq p_m \lt k^{0.45}$ seems more challenging;  we will use results of
Rosser and Schoenfeld \cite{RS} to show this bound for special sums of the above form for $m \geq 184$,
then show how computations bring $m$ down to $20$, and then show how this implies the
general result when $k \gt e^{9.5}$.

We start with getting $\log p_{m+u}$ in terms close to $\log u$. 

\textbf{Lemma}: $\log p_{m+u} \lt \log u  + 
(\log\log u + \frac{m}{u})(1 + \frac{1}{\log u}) + \frac{m}{u(\log u)^2}$ when
$u \gt m$.

Sketch of Proof:
Theorem 6 of \cite{RS} yields for $k \gt 5$ that 
$\log p_k  \lt \log k + \log(\log k + \log\log k)$.
We start with the more complicated subterm:
$\log(\log(m+u) + \log\log(m+u)) = \log\log u + LL$, where
$LL = \log(\log(m+u)/\log u + \log\log(m+u)/\log u)$. 
\begin{eqnarray*}  
  LL &  =  & \log(1 + [\log(1 + m/u) + \log\log(m+u)]/\log u) \\
    &  \lt &  [m/u + \log\log(m+u)]/\log u \\
    &  = & [m/u + \log(\log(u) + \log(1 + m/u))]/\log u \\
    &  = & [m/u + \log\log u + \log(1 + \log(1+ m/u)/\log u )]/\log u \\
    & \lt & [m/u + \log\log u + m/u\log u]/\log u \\
    & = & (m/u + \log\log u)/\log u  + m/u(\log u)^2
\end{eqnarray*}

Now we use Theorem 6 to get
\begin{eqnarray*}
\log p_{m+u} & \lt & \log(m+u) + \log(\log(m+u) + \log\log(m+u)) \\
  & \lt & \log u + m/u + \log\log u + LL \\
  & \lt & \log u + (m/u + \log\log u)(1 + 1/\log u) + m/u(\log u)^2.
\end{eqnarray*}
\textbf{End} of Sketch.

Next, we want to get a good lower bound on $\log(\log p_{m+u}/\log p_m)$ when we have
our special sum close enough to 1.

\textbf{Lemma}: Suppose $\sum_{1 \leq j \leq u} 1/p_{m+j} \gt 1 - 1/2p_m$ and 
$p_{m+u} \gt 286$.  
Then $\log(\log p_{m+u}/\log p_m) \gt 1 - 1/2p_m - .00148 - 1/2(\log p_{m+u})^2$ .

Sketch of Proof: From \cite{RS} Theorem 5 we derive (where
$p \leq x$ means the primes greater than 1 and at most $x$) 
$$\mid \sum_{p\leq x} 1/p - \log\log x -B \mid \leq 1/2(\log x)^2 \textrm{ for all } x \geq 286$$
and from \cite{RS} Theorem 20 we have    
$$ \sum_{p \leq x} 1/p - \log\log x -B \gt 0 \textrm{ for } 2 \leq x \leq 10^8.$$
Combining the results over the two ranges gives  
$$-\log\log p_m - B + 1/2(\log(10^8))^2 >= - \sum_{p \leq p_m} 1/p \textrm{ for all } m \gt 0$$
leading to our estimate:
\begin{eqnarray*}
\sum_{p_m \lt p \leq p_{m+u}} 1/p & \leq & \log\log p_{m+u} + \frac{1}{2(\log p_{m+u})^2} +B \\
                                  &      & -\log\log p_m -B +1/(2(\log 10^8)^2) \\
 & \leq & \log(\log p_{m+u}/\log p_m ) + \frac{1}{2(\log p_{m+u})^2} + .00148
\end{eqnarray*}
Subtracting the last two terms from $1 - 1/2p_m$ leads to the conclusion.
\textbf{End} of Sketch.

Now let us choose $p_m \gt e^7 \gt 1096$ and so $m \geq 184$ and $1/\log u \lt 1/15$.
Then $\log(\log p_{m+u}/\log p_m) \gt 1 - 1/2192 - .00148 - 1/2(15)^2 \gt 0.99$
and $e^{0.99} \gt 2.691$, thus $\log p_{m+u} \gt 2.691 \log p_m$.

Toward a contradiction, assume that $x = \frac{20 \log p_m}{9} \gt \log u$.
Then the Lemma concerning $\log{p_{m+u}}$ gives
\begin{eqnarray*}
1.21095 x & = & 2.691 \log p_m \lt \log p_{m+u} \\
& \lt & x + (\log x + m/u)(1 + 1/\log u) + m/u(\log u)^2 \\
  & \lt & x + (\log x + 1/1000)(1 + 1/15) + 1/1000 \\
  & \lt & x + (16\log x)/15 + 1/100
\end{eqnarray*}
However $.21095 x \gt (16\log x)/15 + .01$ for $x \gt 14$, which means for $\log p_m \gt 7$,
and we have a contradiction.  Thus $u \geq p_m^{20/9}$.

For primes $p_m=71$ through $p_m=857$ ($m=20$ to $m=148$), we verified through 
computation that if
$\sum_{1 \leq j}^u 1/p_{m+i} \gt 1 - 1/2p_m$, then
$\log u /\log p_m \gt 20/9$. In particular for $u$ as low as $13250 \lt e^{9.5}$ and
$n$ a ratio of certain primorials the Improvement holds.  
Also, 
$\sum_{1 \leq j}^{5761308} 1/p_{148+j} \gt 1 - 1/2p_{148}$, thus
$\sum_{1 \leq j}^u 1/p_{m+j} \gt 1 - 1/2p_m$ implies $u \gt p_m^{20/9}$ for $p_m$
running up to $1102 \lt 5761308^{9/20}$.  So the implication holds for all $p_m \geq 71$.  

Now assume $n$ with $k \geq e^{9.5}$ and otherwise arbitrary. Find $l$
 smallest so that $\sum_{l \leq i}^k 1/q_i \lt 1 + 1/2q_l$. Let $p_m=q_{l-1}$. Then 
$1 - 1/2q_{l-1} \leq \sum_{l \leq i}^k 1/q_i \leq 
\sum_{1 \leq j}^{k-l+1} 1/p_{m+j}.$
Either $p_m \geq 71$ and so $(k-l+1)^{0.45}\geq u^{0.45} \geq p_m$, or else 
$q_{l-1} \lt 71 \lt (e^{9.5})^{9/20} \leq k^{0.45}$.

\section{Addendum}

Shortly after version 1 of this article was posted, we found an
upper bound based on estimates of Euler's totient which does 
better than the Improvement.

Aaron Meyerowitz asked about these estimates in question 88777 
on MathOverflow.  The
key observation is that the number of totatives to $n$ in the
interval $[0,x)$ differs from $x\phi(n)$ by a periodic function 
$E_n(x)$ which has a maximum value at most $2^{k - 1}$.  Using 
this, we can generalize the Variation using a similar argument 
to show
$$ g(n) \leq 
\frac{(k-m+1)(2^m + \pi^{-1}(d))}{\pi^{-1}(d)(1 - \sigma^{-1}(n/d))},$$
where one chooses $m \lt k$ so that $d=\prod_{1 \leq i \leq m} q_i$ 
divides $n$ and also satisfies
$t=1-\sigma^{-1}(n/d) \gt 0$, and finally so that the 
right hand side above leads to an optimal bound. (Hint:
the Proof sketch of the Variation introduces an error term of 
size at most $2$ for counting multiples of $q_i$ which are not 
multiples of $q_1$; now use an error of $2^m + \pi^{-1}(d)$ 
for an interval $[y,x]$ of length $L$ when counting multiples 
of $q_i$ that are coprime to $d$.)

When such an optimal $m>0$ is found, one can easily show that 
$t > 1/(3+q_{m+1})$, and the analysis from the Appendix can 
be adapted to show $m^e \lt k$.  As $1/\pi^{-1}(d)$ is bounded by 
$2\log q_m$, this will beat our bound $k^{K\sigin}$ when 
$2^m$ is substantially smaller than $k^{K\sigin -3/2}$. 
Computing both bounds for $n=P_k$ for $k<10^7$ show this new
bound to be superior: we expect to show that it holds when
$k < 10^{10}$ and $\sigin > 1$, which would imply that $K$ above
can be taken near $3.6$ for all $k \gt 2$.

This bound can be used to show Kanold's bound $2^{\sqrt{k}}$ 
actually holds for $k \gt e^6$.  There are also improvements
to be made on the error term $(2^m + \pi^{-1}(d))$; with
such improvements we expect to show in a followup article a 
subquadratic in $k$ upper bound for $k< 10^{10}$.

We also found a statement of the Observation in a 1975 work of
Kanold's.  (We thank Prof. Dr. Heiko Harborth for making this
part of the literature available to us.) 
We are still looking for a published proof of the
Observation as well as an appearance in the literature of the
Proposition that $g(n) \leq \myceil{k/(1-\sigin)}$ for those $n$ with
$\sigin \lt 1$.  We still believe that the explicit upper bound of
$k^{3 + 3.81\sigin}$ has not appeared in the literature.

We have planned a series of forthcoming articles, tentatively titled
'Adventures in finding bounds on Jacobsthal's function.'  In addition
to fleshing out some of the questions asked in an earlier section, we
will consider the computational complexity of $g(n)$ and approximations to
$g(n)$,
various lower bounds coming from elementary (and not so elementary)
considerations, applications using both conjectured and actual bounds,
and generalizations in geometric and algebraic realms.

\end{document}